\documentclass{article}
\usepackage[latin1]{inputenc}
\usepackage{amsmath,amsthm,amssymb}

\newtheorem{thm}{Theorem}
\newtheorem{lemma}[thm]{Lemma}

\def\eq#1{(\ref{#1})}

\newcommand{\neweq}[1]{\begin{equation}\label{#1}}

\def\ep{\varepsilon}

\def\phi{\varphi}
\def\RR{\mathbb R}

\def\di{\displaystyle}
\def\ri{\rightarrow}

\def\incep{\left\{\begin{array}{cl} }
 \def\termin{\end{array}\right. }
\def\2af{2^*_\alpha}

\title{Singular solutions of perturbed logistic-type
equations}
\author{\sc Du\v{s}an Repov\v{s}$\,^{a,b}$\\
 \small $^a\,$Faculty of Mathematics and
Physics, University of Ljubljana,\\ \small Jadranska  19,  P. O. Box 2964, 1001 Ljubljana, Slovenia\\
\small $^b\,$Faculty of Education, University of Ljubljana,\\ \small Kardeljeva plo\v{s}\v{c}ad 16, 1000 Ljubljana, Slovenia\\
 \small E-mail: {\tt dusan.repovs@guest.arnes.si}\\}

\begin{document}

\maketitle
\begin{abstract}
We are concerned with the qualitative analysis of positive singular solutions with blow-up boundary for a class of logistic-type equations with slow diffusion and variable potential. We establish the exact blow-up rate of solutions near the boundary in terms of Karamata regular variation theory. This enables us to deduce the uniqueness of the singular solution.
 \\
{\bf Keywords:} singular solution, blow-up boundary, logistic equation, Karamata regular variation theory, maximum
principle.\\
{\bf 2010 Mathematics Subject Classification:} 35J60; 35B40; 35B50.
\end{abstract}

\section{Introduction}
Let $\Omega$ be a bounded domain with smooth boundary in $\RR^N$, $N\geq 1$. Assume $f:[0,\infty)\ri [0,\infty)$ is a locally Lipschitz continuous function such that
\begin{equation}\label{h1}
f(0)=0\qquad\mbox{and}\qquad f(t)>0\ \ \mbox{for}\ \ t>0\end{equation}
and
\begin{equation}\label{h2}
f\mbox{ is nondecreasing}.\end{equation}

Consider the basic population model described by the logistic problem
\begin{equation}\label{p0}\left\{
\begin{array}{lll}
&\Delta u=f(u)&\qquad\mbox{in $\Omega$}\\
&\lim_{x\ri\partial\Omega}u(x)=+\infty\\
&u>0&\qquad\mbox{in $\Omega$}\end{array}\right.
\end{equation}

All smooth functions satisfying problem \eq{p0} are called {\it large} (or {\it blow-up boundary}) solutions.

Under assumptions \eq{h1} and \eq{h2}, Keller \cite{keller} and Osserman \cite{osserman} proved that problem \eq{p0} has a solution if and only if
\begin{equation}\label{h3}\int^{+\infty}\frac{1}{\sqrt{F(u)}}\,du<+\infty ,\end{equation}
where $F(u):=\int_0^uf(s)ds$.

We refer to Ghergu and R\u adulescu \cite[Theorem 1.1]{grcpaa} for an elementary argument that problem \eq{p0} cannot have any solution if $f$ has a sublinear or a linear growth, hence it does not satisfy condition \eq{h3}.
We point out that the original approach is due to Dumont, Dupaigne, Goubet and R\u adulescu \cite{ddgr}, who removed the monotonicity assumption \eq{h2} and showed that the key role in the existence of solutions of problem \eq{p0} is played only by the {\it Keller-Osserman condition} \eq{h3}.

Functions satisfying the Keller-Osserman condition have a super-linear growth, such as: (i) $f(u)=u^p$ ($p>1$); (ii) $f(u)=e^u$; (iii) $f(u)=u^p\ln (1+u)$ ($p>1$); (iv) $f(u)=u\ln^p(1+u)$ ($p>2$).

We point out that the study of large solutions was initiated by Bieberbach \cite{bieb} in 1916 and Rademacher \cite{rade} in 1943 for the particular case $f(u)=e^u$ if $N=2$ or $N=3$.
An important contribution to the study of singular solutions with boundary blow-up is due to Loewner and Nirenberg \cite{loni}, who linked the uniqueness of the large
solution to the growth rate at the boundary. Motivated by certain geometric problems, they
established the uniqueness of the solution in the case $f(u)=u^{(N+2)/(N-2)}$, $N\geq 3$.

 C\^{\i}rstea and  R\u adulescu studied in \cite{crhouston} (see Du and Guo \cite{duguo} for the quasilinear case) the perturbed logistic problem
\begin{equation}\label{p1}\left\{
\begin{array}{lll}
&\Delta u+au=b(x)f(u)&\qquad\mbox{in $\Omega$}\\
&\lim_{x\ri\partial\Omega}u(x)=+\infty\\
&u>0&\qquad\mbox{in $\Omega$},\end{array}\right.
\end{equation}
where $a$ is a real number and $b\in C^{0,\alpha}(\overline{\Omega})$, $0<\alpha<1$, such that $b\geq 0$ and $b\not\equiv 0$ in $\Omega$. C\^{\i}rstea and  R\u adulescu found the whole range of values of the
parameter $a$ such that problem \eq{p1} admits a solution and this responds to a question raised by H.~Brezis. Their analysis includes the case where the potential $b(x)$
vanishes on $\partial\Omega$. Due to the fact that $u$ has a singular behavior on the boundary, this setting corresponds to the ``competition" $0\cdot\infty$ on $\partial\Omega$. The study carried out in \cite{crhouston} strongly relies on the structure of the subset of $\Omega$ where the potential $b$ vanishes. In particular, it is argued in \cite{crhouston} that problem \eq{p1} has a solution for all values of $a\in\RR$ provided that
$$\mbox{int}\,\{x\in\Omega;\ b(x)=0\}=\emptyset.$$
We also refer to Ghergu and R\u adulescu \cite{groxford} for related results.

Our main purpose in this paper is to study the effect of a {\it sublinear} perturbation $au^p$ ($0<p<1$) in problem \eq{p0}. This framework corresponds to a {\it slow diffusion} in the population model.
According to Delgado and Su\'arez, the assumption $0<p<1$ means that the
diffusion, namely the rate of movement of the species from high density regions to low density ones, is
slower than in the linear case corresponding to $p=1$, which is described by problem \eq{p1}.

\section{Statement of the problem and main results}
We start with the following example of singular logistic indefinite superlinear model.
Fix $m>1$ and consider the nonlinear problem
\begin{equation}\label{pmodel}\left\{
\begin{array}{lll}
&\Delta w^m+aw=b(x)w^2&\qquad\mbox{in $\Omega$}\\
&\lim_{x\ri\partial\Omega}w(x)=+\infty\\
&w>0&\qquad\mbox{in $\Omega$}.\end{array}\right.
\end{equation}

This problem can be regarded as a model of a steady-state single species inhabiting in
$\Omega$, so $w(x)$ stands for the population density. The parameter $a$ represents the growth rate
of the species while
the term $m > 1$ was introduced by Gurtin and  MacCamy \cite{gurtin} to describe the dynamics of biological
population whose mobility depends upon their density. We refer to Li,  Pang and Wang \cite{lipang} for a study of problem \eq{pmodel}
in the case of multiply connected domains and subject to mixed boundary conditions.

The change of variable $u=w^m$ transforms problem \eq{pmodel} into
\begin{equation}\label{pmodel1}\left\{\begin{array}{lll}
&\Delta u+au^p=b(x)u^q&\qquad\mbox{in $\Omega$}\\
&\lim_{x\ri\partial\Omega}u(x)=+\infty\\
&u>0&\qquad\mbox{in $\Omega$},\end{array}\right.\end{equation}
where $p=1/m\in (0,1)$ and $q=2/m$. As stated in the previous section, it is expected that this problem has a solution in the super-linear setting, that is, provided that $m<2$.

In this paper we study the more general problem
$$\left\{\begin{array}{lll}
&\Delta u+ag(u)=b(x)f(u)&\qquad\mbox{in $\Omega$}\\
&\lim_{x\ri\partial\Omega}u(x)=+\infty\\
&u>0&\qquad\mbox{in $\Omega$},\end{array}\right.$$
where $g$ has a sublinear growth and $f$ is a function satisfying the Keller-Osserman condition such that the mapping $f/g$ is increasing in $(0,\infty)$. To fix the ideas, we consider the
model problem
\begin{equation}\label{p}\left\{\begin{array}{lll}
&\Delta u+au^p=b(x)f(u)&\qquad\mbox{in $\Omega$}\\
&\lim_{x\ri\partial\Omega}u(x)=+\infty\\
&u>0&\qquad\mbox{in $\Omega$}.\end{array}\right.\end{equation}

In order to describe our main result we recall some basic notions and properties of the Karamata theory of functions with regular variation at infinity.
We refer to Bingham, Goldie and Teugels \cite{bingham} and Seneta \cite{seneta} for more details.

A positive
measurable function $R$ defined on $[A,\infty)$, for
some
$A>0$, is called {\it regularly varying} (at infinity)
with index  $q\in \RR$, written $R\in \RR_q$, if
for all $\xi>0$
$$ \lim_{u\to \infty}R(\xi u)/R(u)=\xi^q. $$
If $R:[A,\infty)\to
(0,\infty)$
is measurable and Lebesgue integrable on each finite
subinterval of
$[D,\infty)$, then $R$ varies regularly if and only if
there exists $j\in \RR$ such that
\begin{equation}\label{lt} \lim_{u\to
\infty}\frac{u^{j+1}R(u)}{\int_D^u x^j R(x)\,dx}
\end{equation}
exists and is a positive number, say $a_j+1$.
In this case, $R\in \RR_q$ with $q=a_j-j$.
Moreover, by a theorem established in 1933 by Karamata,
if $R\in \RR_q$ is Lebesgue integrable
on each finite subinterval of $[D,\infty)$, then the
limit defined by
\eq{lt} is $q+j+1$, for every $j>-q-1$.
We also point out that if $S\in C^1[A,\infty)$, then
$S'\in \RR_q$ with $q>-1$ if and only if, for some $m>0$, $C>0$
and $B>A$, we have
$$ S(u)= Cu^m {\rm exp}\left\{ \int_B^u\frac{y(t)}{t}\,dt\right\}\qquad\mbox{for all}\
 u\geq B ,$$ where $y\in C[B,\infty)$ satisfies $\lim_{u\to \infty}
y(u)=0$. In this case, $S'\in \RR_q$ with $q=m-1$.

As established in C\^{\i}rstea and R\u adulescu \cite{crcras}, if
$ f'\in \RR_\rho$ then $\rho\geq 0$ and, furthermore, if $\rho>0$ then $f$ satisfies the Keller-Osserman condition, provided that $f$ is increasing.

Next, we denote by
${\cal K}$ the Karamata class containing all positive, increasing
$C^1$-functions $k$
defined on $(0,\nu)$, for some $\nu>0$, which satisfy
$\lim_{t\to 0^+}\di\left(\frac{\int_0^t
k(s)\,ds}{k(t)}\right)^{(i)}:=\ell_i,\
i=\overline{0,1}$.
A straightforward computation shows that
$\ell_0=0$ and $\ell_1\in [0,1]$, for every $k\in {\cal K}$.
We refer to Lemma 2 in \cite{crcras}, where it is argued that $\ell_1$ can actually take any value in $[0,1]$.

Throughout this work we assume that $a$ is a real parameter and $b\in C^{0,\alpha}(\overline{\Omega})$, $0<\alpha<1$, such that $b\geq 0$ and $b\not\equiv 0$ in $\Omega$. We also assume that $f:[0,\infty)\ri [0,\infty)$ is a locally Lipschitz continuous function that satisfies hypotheses \eq{h1}, \eq{h3} and
\begin{equation}\label{h4}
\mbox{the mapping}\ (0,\infty)\ni u\longmapsto\frac{f(u)}{u^p}\ \mbox{is increasing}.\end{equation}

Our first result establishes the existence of a unique positive singular solution of problem \eq{p}. The existence is deduced by means of a suitable comparison principle in the case of semilinear elliptic equations without boundary condition.  Next, we use this existence result to study the same nonlinear elliptic equation with sublinear perturbation in the framework of non simply connected domains and subject to mixed boundary condition. In both cases, the uniqueness of the solution follows after establishing the blow-up rate of an arbitrary solution  near the boundary. Throughout this paper we denote $d(x):=\mbox{dist}\, (x,\partial\Omega)$, for all $x\in\Omega$.

\begin{thm}\label{t1}
Assume conditions \eq{h1}, \eq{h3} and \eq{h4} are fulfilled. Then problem \eq{p} has at least one solution.

Assume hypotheses \eq{h1}, \eq{h4} and  $ f\in \RR_{\rho+1}$ with $\rho>0$. Assume the potential $b(x)$ satisfies
$$b(x)=c\,k^2(d(x))+o(k^2(d(x))\qquad\mbox{as}\ d(x)\ri 0,$$
where $c$ is a positive number and $k\in {\mathcal K}$. Then, for all real number $a$, problem \eq{p} has a unique solution $u_a$ and
\begin{equation}\label{58}
u_a(x)=\xi_0\,h(d(x))+o(h(d(x))\qquad\mbox{as $d(x)\ri 0$},\end{equation}
where $h$ is uniquely defined by
\begin{equation}\label{63}\int_{h(t)}^\infty\frac{ds}{\sqrt{F(s)}}=\sqrt 2\int_0^tk(s)ds\end{equation}
and
$$\xi_0=\left(\frac{2+\rho\ell_1}{(2+\rho )c}\right)^{1/\rho}.$$
\end{thm}

The existence result described in the first part of Theorem \ref{t1} is in contrast with the corresponding one for the linear perturbed case studied in C\^{\i}rstea and R\u adulescu \cite{crhouston}. In their analysis a key role is played by the set
$\Omega_0:=\mbox{int}\,\{x\in\Omega;\ b(x)=0\}.$
Let $H_{\infty}$ define
the Dirichlet Laplacian on the set $\Omega_0\subset\subset
\Omega$
as the unique self-adjoint operator associated to the
quadratic
form $\psi(u)=\int_{\Omega}|\nabla u|^2\,dx$
with form domain
$$ H^1_D(\Omega_0)=\{u\in H_0^1(\Omega);\ u(x)=0\quad
\mbox{for a.e.}\ x\in \Omega\setminus\Omega_0\}. $$
If $\partial\Omega_0$ satisfies an exterior
cone condition, then  $H^1_D(\Omega_0)$
coincides with $H^{1}_{0}(\Omega_0)$ and $H_{\infty}$
is the classical Laplace operator with Dirichlet
condition on $\partial\Omega_0$.
Let $\lambda_{\infty,1}$ be
the first Dirichlet eigenvalue of  $H_\infty$ in
$\Omega_0$. If $\Omega_0=\emptyset$ then
 $\lambda_{\infty,1}=+\infty$.

The main result in \cite{crhouston} asserts that problem \eq{p} (for $p=1$) has a solution if and only if $a\in
(-\infty,\lambda_{\infty,1})$. By contrast, our result established in Theorem \ref{t1} shows that the perturbation $au^p$ is {\it small} enough provided that $0<p<1$ and this does not affect the existence of blow-up boundary solutions in the sublinear setting
we have described above.  The basic assumption $0<p<1$ allows us to construct an appropriate super-solution for problem \eq{p} for all $a\in\RR$. The similar argument is no more possible provided that $p=1$.

We can also see problem \eq{p} as a perturbation of the blow-up boundary logistic equation
$$\left\{\begin{array}{lll}
&\Delta u=b(x)f(u)&\qquad\mbox{in $\Omega$}\\
&\lim_{x\ri\partial\Omega}u(x)=+\infty\\
&u>0&\qquad\mbox{in $\Omega$},\end{array}\right.$$
where $f$ is a positive increasing function satisfying the Keller-Osserman condition. Combining our result with those obtained in C\^{\i}rstea and R\u adulescu \cite{crcras,crhouston} we may assert the following: (i) in the sublinear case  $0<p<1$, the perturbed equation \eq{p} has a unique solution for all $a\in\RR$; (ii) in the linear perturbed case $p=1$, the problem \eq{p} has a (unique) solution if and only if $a<\lambda_{\infty,1}\leq +\infty$.

Next, we assume that $\emptyset\not=\overline{\Omega_0}\subset\Omega$. We denote $D:=\Omega\setminus \overline{\Omega_0}$
and we assume that $b>0$ in $D$. We are now concerned with the nonlinear problem
\begin{equation}\label{pp}\left\{\begin{array}{lll}
&\Delta u+au^p=b(x)f(u)&\qquad\mbox{in $D$}\\
&u=0&\qquad\mbox{on $\partial\Omega$}\\
&\lim_{x\ri\partial\Omega_0}u(x)=+\infty\\
&u>0&\qquad\mbox{in $D$}.\end{array}\right.\end{equation}

It is striking to remark that solutions of problem \eq{pp} fulfill similar properties as those established in Theorem \ref{t1}. A related result may be found in C\^{\i}rstea and R\u adulescu \cite{crcras2}.

\begin{thm}\label{t2}
Assume conditions \eq{h1}, \eq{h3} and \eq{h4} are fulfilled. Then problem \eq{pp} has a minimal and a maximal solution.

Assume hypotheses \eq{h1}, \eq{h4} and  $ f\in \RR_{\rho+1}$ with $\rho>0$. Assume the potential $b(x)$ satisfies
$$b(x)=c\,k^2(d(x))+o(k^2(d(x))\qquad\mbox{as}\ d(x)\ri 0,$$
where $c$ is a positive number and $k\in {\mathcal K}$. Then, for all real number $a$, problem \eq{pp} has a unique solution $u_a$ and
$$
u_a(x)=\xi_0\,h(d(x))+o(h(d(x))\qquad\mbox{as $d(x)\ri 0$},$$
where $h$ is uniquely defined by
$$\int_{h(t)}^\infty\frac{ds}{\sqrt{F(s)}}=\sqrt 2\int_0^tk(s)ds$$
and
$$\xi_0=\left(\frac{2+\rho\ell_1}{(2+\rho )c}\right)^{1/\rho}.$$
\end{thm}

Theorems \ref{t1} and \ref{t2} can be extended to a Riemannian manifold setting if the Laplace operator is replaced by the Laplace-Beltrami differential operator
$$\Delta_B:=\frac{1}{\sqrt d}\,\frac{\partial}{\partial x_i}\left(\sqrt d\, a_{ij}(x)\,
\frac{\partial}{\partial x_i}\right),\qquad d:=\mbox{det}(a_{ij})$$
with respect to the metric $ds^2=b_{ij}dx_idx_j$, where $(b_{ij})$ denotes the inverse of $(a_{ij})$. We refer, e.g., to Loewner and Nirenberg \cite{loni}, where $\Omega$ is replaced by the sphere $(S^N,g_0)$ and $\Delta$ is the Laplace-Beltrami operator $\Delta_{g_0}$.

\section{Proofs of the main results}
A central role is played by the following comparison principle for logistic-type equations with sublinear perturbation. The proof relies on some ideas introduced by Benguria, Brezis and Lieb \cite{bbl} (see also Marcus and V\'eron
\cite[Lemma 1.1]{marver}, C\^{\i}rstea and R\u adulescu \cite[Lemma 1]{crhouston}, Du and Guo \cite{duguo}).

\begin{lemma}\label{l1}
Let $D$ be a bounded domain in $\RR^N$ with smooth boundary. Assume $a$ is a real number and let $h$, $r$ be $C^{0,\alpha}$-functions in $\overline D$ such that $h\geq 0$ and $r\geq 0$ in $D$. Let $u_1,\, u_2\in H^1(D)$ be positive continuous functions such that
\begin{equation}\label{l11}\begin{array}{ll}
&\di \Delta u_1+au_1^p-h(x)f(u_1)+r(x)\leq 0\leq \\ &\di\Delta u_2+au_2^p-h(x)f(u_2)+r(x)\qquad\mbox{in ${\mathcal D}'(D)$}\end{array}
\end{equation}
and
\begin{equation}\label{l12}
\limsup_{x\ri\partial D}(u_2(x)-u_1(x))\leq 0,
\end{equation}
where $f$ is continuous on $[0,\infty)$ such that the mapping $f(t)/t^p$ is increasing for $\inf_D(u_1,u_2) <t<
\sup_D(u_1,u_2)$.

Then $u_1\geq u_2$ in $D$.
\end{lemma}

\proof
Relation \eq{l11} implies that for all $\psi\in C^2_c(D)$ with $\psi\geq 0$ we have
\begin{equation}\label{l13}\begin{array}{ll}
&\di\int_D\left(\nabla u_1\nabla\psi-au_1^p\psi +hf(u_1)\psi -r\psi\right)dx \geq 0\geq\\ &\di
\int_D\left(\nabla u_2\nabla\psi-au_2^p\psi +hf(u_2)\psi -r\psi\right)dx.\end{array}
\end{equation}

Relation $u_1\geq u_2$ in $D$ is equivalent to $G:=\{x\in D;\ u_1(x)<u_2(x)\}=\emptyset$.
Fix $\ep>0$ small enough and denote
$$D(\ep):=\{x\in D;\ u_2(x)>u_1(x)+\ep\}.$$
For $i=1,2$ we set
$$v_i=(u_i+\ep_i)^{-p}\left( (u_2+\ep_2)^{1+p}-(u_1+\ep_1)^{1+p}\right)^+,$$
where $\ep_1=2\ep$, $\ep_2=\ep$.
Thus, $v_i\in H^1(D)$ and it vanishes outside the set $D$. Using ow assumption  \eq{l12}, we have $D(\ep)\subset\subset D$. Hence, $v_i$ can be approximately in the $H^1\cap L^\infty$ topology on $\overline D$ by nonnegative $C^2$ functions vanishing near $\partial D$. It follows that relation \eq{l13} holds with $v_i$ instead of $\psi$.
We deduce that
\begin{equation}\label{l14}\begin{array}{ll}
&\di\int_{D(\ep)}\left(\nabla u_2\nabla v_2-\nabla u_1\nabla v_1\right)dx +
\int_{D(\ep)}h(x)\left( f(u_2)v_2-f(u_1)v_1\right)dx\leq\\
&\di  \int_{D(\ep)} a(u_2^pv_2-u_1^pv_1)dx+\int_{D(\ep)} r(x)(v_2-v_1)dx.\end{array}
\end{equation}
With a straightforward computation, as in the proof of Lemma 1 in \cite{crhouston}, we deduce that
\begin{equation}\label{l144}\begin{array}{ll}
&\di \nabla u_2\nabla v_2-\nabla u_1\nabla v_1=\left[1+p\left(\frac{u_2+\ep}{u_1+2\ep}\right)^{1+p}\right]|\nabla u_1|^2+\\ &\di \left[1+p\left(\frac{u_1+2\ep}{u_2+\ep}\right)^{1+p}\right]|\nabla u_2|^2-\\ &\di
(1+p)\left[\left(\frac{u_2+\ep}{u_1+2\ep}\right)^{p}+\left(\frac{u_1+2\ep}{u_2+\ep}\right)^{p}\right]\nabla u_1\cdot\nabla u_2\geq 0.\end{array}
\end{equation}
Since $f(t)/(t+\ep)^p$ is increasing on $(0,\infty)$, we find
$$\frac{f(u_1)}{(u_1+2\ep)^p}<\frac{f(u_1+\ep)}{(u_1+2\ep)^p}<\frac{f(u_2)}{(u_2+\ep)^p}\qquad\mbox{in $D(\ep)$}.$$
Thus, all the integrands on the left-hand side of \eq{l14} are nonnegative,
while the second term on the right-hand side of \eq{l14} equals to
\begin{equation}\label{l145}-\int_{D(ep)}r(x)\,\frac{\left[(u_2+\ep)^{1+p}-(u_1+2\ep)^{1+p}\right]\left[(u_2+\ep)^{p}-
(u_1+2\ep)^{p}\right]}{(u_1+2\ep)^{p}(u_2+\ep)^{p}}\,dx\leq 0.\end{equation}

Relations \eq{l14}, \eq{l144} and \eq{l145} show that $\limsup_{\ep\ri 0}A_\ep\geq 0$, where
$$A_\ep:= \int_{D(\ep)}(u_2^pv_2-u_1^pv_1)dx.$$
A precise answer is given in what follows. We point out that the result stated below is obvious in the linear case that corresponds to $p=1$, see C\^{\i}rstea and  R\u adulescu \cite{crhouston}.

{\sc Claim}. {\it We have $\lim_{\ep\ri 0}A_\ep=0$.}

{\it Proof of Claim}. Fix $\eta>0$ and $\rho>0$ such that $\rho^p(1+\rho)<\eta$. Set $$D_1(\ep,\rho):=\{x\in D(\ep);\ u_2(x)<\rho\}\qquad \mbox{and}\qquad D_2(\ep,\rho):=D(\ep)\setminus D_1(\ep,\rho).$$
We first observe that for all $\ep\in(0,1)$,
$$\begin{array}{ll}&\di\int_{D_1(\ep,\rho)} (u_2^pv_2-u_1^pv_1)dx=\\&\di
\int_{D_1(\ep,\rho)} \left(\frac{u_2^p}{(u_2+\ep)^p}-\frac{u_1^p}{(u_1+2\ep)^p} \right)\left[(u_2+\ep)^{1+p}-(u_1+2\ep)^{1+p})\right]dx\leq\\
&\di  \int_{D_1(\ep,\rho)} \left[ u_2^p(u_2+\ep)+u_1^p(u_1+2\ep)\right]dx\leq C(\Omega,k)\eta.\end{array}$$

In order to estimate the integral of the same quantity over $D_2(\ep,\rho)$ we first observe that there is some $C_1>0$ such that for all $\ep\in (0,1)$ and for any $x\in D_2(\ep,\rho)$, we have
$u_1^p(x)/(u_1(x)+2\ep)\leq C_1$. This follows by a contradiction argument combined with the assumptions $0<p<1$ and $u_2\geq\rho$ on $D_2(\ep,\rho)$. On the other hand, with the same arguments as in Du and Guo \cite[p.~283]{duguo}, there exists $C_2>0$ such that for all $\ep\in (0,1)$ and for any $x\in D_2(\ep,\rho)$,
$$(u_2+\ep)^{1+p}-(u_1+2\ep)^{1+p}\leq C_2.$$
This enables us to conclude that
$$\limsup_{\ep\ri 0}\int_{D_2(\ep,\rho)}\left(\frac{u_2^p}{(u_2+\ep)^p}-\frac{u_1^p}{(u_1+2\ep)^p} \right)\left[(u_2+\ep)^{1+p}-(u_1+2\ep)^{1+p})\right]dx\leq 0.$$

Therefore
$$0\leq \liminf_{\ep\ri 0}A_\ep\leq \limsup_{\ep\ri 0}\int_{D_1(\ep,\rho)} (u_2^pv_2-u_1^pv_1)dx\leq C(\Omega,k)\eta,$$
for all $\eta>0$.
Since $\limsup_{\ep\ri 0}A_\ep\geq 0$, our claim follows.

\smallskip
It now remains to observe that the set $G$ is empty.
We argue by contradiction and assume that $G\not=\emptyset$. Fix arbitrarily $x_0\in G$ and take a small closed ball $B$ centered at $x_0$ such that $B\subset G$. Since $\min_B(u_2-u_1)=:m>0$, we deduce that $B\subset D(\ep)$ for all $\ep\in (0,m)$. But
$$\begin{array}{ll}
0&\di\leq\int_{B}\left(\nabla u_2\nabla v_2-\nabla u_1\nabla v_1\right)dx +
\int_{B}h(x)\left( f(u_2)v_2-f(u_1)v_1\right)dx\\
&\di -\int_{B} r(x)(v_2-v_1)dx\leq a\int_{D(\ep)} (u_2^pv_2-u_1^pv_1)dx.\end{array}$$
Letting $\ep\ri 0^+$ we deduce that for all $x\in B$,
$$\frac{\nabla u_1(x)}{u_1(x)}=\frac{\nabla u_2(x)}{u_2(x)}\qquad\mbox{and}\qquad h(x)=0.$$
Since $x_0\in G$ is arbitrary, we obtain $\nabla (\ln u_2-\ln u_1)=0$ and $h\equiv 0$ in $G$. But $h\not\equiv 0$ in $D$, hence $G\not=D$. Thus, $\partial G\cap D\not=\emptyset$. We take $x'\in \partial G\cap D$ and $\omega\subset G$ such that $x'\in\partial\omega$. Hence $u_1(x')=u_2(x')$ and $\nabla (\ln u_2-\ln u_1)\equiv 0$ in $\omega$, hence $u_2/u_1\equiv C>0$ in $\omega$. By continuity we deduce that $C=1$, which shows that $u_1=u_2$ in $\omega$. This contradicts $\omega\subset G$. Thus, $u_1\geq u_2$ in $D$ and this concludes the proof of our lemma.
\qed

\begin{lemma}\label{l2}
Let $D$ be a bounded domain in $\RR^N$ with smooth boundary. Assume $h$, $k$ and $r$ are $C^{0,\alpha}$-functions in $\overline D$ such that $h>0$, $k\geq 0$ and $r\geq 0$ in $D$. Then for any non-negative function $0\not\equiv \Phi\in
C^{0,\alpha} (\partial D)$, the nonlinear problem
\begin{equation}\label{l21}
\left\{
\begin{array}{lll}
&\di \Delta u+k(x)u^p=h(x)f(u)-r(x)&\qquad\mbox{in $D$}\\
&\di u>0&\qquad\mbox{in $D$}\\
&\di u=\Phi &\qquad\mbox{on $\partial D$}
\end{array}\right.
\end{equation}
has a unique solution.
\end{lemma}

\proof We first observe that, by Lemma \ref{l1}, problem \eq{l21} has at most one solution. To prove the existence of a solution we use the method of lower and upper solutions. Due to the sublinear perturbation $k(x)u^p$, the construction provided in the proof of Lemma 2 in \cite{crhouston} does not apply to our framework. However, we observe that $\underline U(x)=0$ is a sub-solution of \eq{l21}. Next, we construct a positive super-solution of \eq{l21} and, by the maximum principle, we argue that this solution is positive in $D$.

Consider $\overline U(x)=M\varphi_1(x)$, where $M>0$ is big enough and $\varphi_1>0$ is an eigenfunction of the Laplace operator in $H^1_0(\Omega)$, where $\Omega\subset\RR^N$ is a smooth bounded domain such that $D\subset\subset\Omega$.
Since $0<p<1$ and $M>0$ is large, we deduce that $\overline U$ is a super-solution of problem \eq{l21}. Thus, problem \eq{l21} has a solution $u_0$ such that $0\leq u_0\leq M\varphi_1$ in $D$.
By standard Schauder and H\"older bootstrap arguments, $u_0$ is a classical solution of problem \eq{l21}.

It remains to argue that $u_0>0$ in $D$. Indeed, since $u_0\leq M\varphi_1$ and the
 mapping $f(u)/u^p$ is increasing on $(0,\infty)$, there is some $C_0>0$ such that $h(x)f(u_0)\leq C_0$ for all $x\in\Omega$. Therefore
$$-\Delta u_0+C_0\geq -\Delta u_0+h(x)f(u_0)=h(x)u_0^p+r(x)\geq 0\qquad\mbox{in $D$}.$$
Since $u_0\geq 0$ on $\partial D$, the maximum principle (see Pucci and Serrin \cite{pucser}) implies that $u_0>0$ in $D$. This completes the proof of Lemma \ref{l2}.
\qed

By taking $\Phi (x)=n$ in \eq{l21} we obtain a sequence of corresponding solutions $(u_n)$ such that $u_n\leq u_{n+1}$ in $D$. We now argue that $(u_n)$ is locally bounded in $D$ provided that $f$ satisfies hypotheses  \eq{h1}, \eq{h3} and \eq{h4}. Indeed, let $\overline u$ be a solution of the singular problem
\begin{equation}\label{amiens}\left\{
\begin{array}{lll}
&\Delta u=\underline hf(u)-\overline k-\overline r-1&\qquad\mbox{in $D$}\\
&\lim_{x\ri\partial D}u(x)=+\infty\\
&u>0&\qquad\mbox{in $D$},\end{array}\right.
\end{equation}
where $\underline h=\min_{\overline D}h(x)$, $\overline k=\max_{\overline D}k(x)$, and $\overline r=\max_{\overline D}r(x)$.
Such a solution exists according to the general results established in Dumont, Dupaigne, Goubet and R\u adulescu \cite{ddgr}.
By the maximum principle, $u_n\leq u_{n+1}\leq \overline u$ in $D$. Thus, under the assumptions of Lemma \ref{l2} and if satisfies hypotheses  \eq{h1}, \eq{h3} and \eq{h4}, we deduce that $(u_n)$ converges to a solution of the singular problem
$$\left\{
\begin{array}{lll}
&\Delta u+k(x)u^p=h(x)f(u)-r(x)&\qquad\mbox{in $D$}\\
&\lim_{x\ri\partial D}u(x)=+\infty\\
&u>0&\qquad\mbox{in $D$}.\end{array}\right.
$$

\begin{lemma}\label{l3}
Let $\Omega$ be a bounded domain in $\RR^N$ with smooth boundary. Assume conditions \eq{h1}, \eq{h3} and \eq{h4} are fulfilled. Let $0\not\equiv \Phi\in
C^{0,\alpha} (\partial \Omega)$ be a  non-negative function  and
$b\in
C^{0,\alpha} ( \Omega)$ be such that $b\geq 0$ in $\Omega$ and $b>0$ on $\partial\Omega$. Then
the nonlinear problem
\begin{equation}\label{l31}
\left\{
\begin{array}{lll}
&\di \Delta u+au^p=b(x)f(u)&\qquad\mbox{in $\Omega$}\\
&\di u>0&\qquad\mbox{in $\Omega$}\\
&\di u=\Phi &\qquad\mbox{on $\partial \Omega$}
\end{array}\right.
\end{equation}
has a unique solution for all $a\in\RR$ and $0<p<1$.
\end{lemma}

\proof We follow an idea developed in the proof of Lemma 3 in C\^{\i}rstea and R\u adulescu \cite{crhouston}.

We first observe that, by Lemma \ref{l1}, problem \eq{l31} has at most one solution.

{\sc Case 1}: $a\geq 0$.

We first observe that the function $\underline U=0$ is a lower solution of problem \eq{l31}.

Let $\Omega_i$ ($i=0,1,2$) be sub-domains of $\Omega$ with smooth boundaries such that $\Omega_0\subset\subset\Omega_1\subset\subset\Omega_2\subset\subset\Omega$. The above remark shows that the nonlinear singular problem
$$
\left\{
\begin{array}{lll}
&\di \Delta u+au^p=b(x)f(u)&\qquad\mbox{in $\Omega\setminus\overline{\Omega_1}$}\\
&\di u>0&\qquad\mbox{in $\Omega\setminus\overline{\Omega_1}$}\\
&\di u=+\infty &\qquad\mbox{on $\partial \Omega\cup\partial\Omega_1$}
\end{array}\right.
$$
has a solution $u_\infty$. Next, we construct a function $u_+\in C^2(\Omega)$ such that $u_+=u_\infty$ in
$\Omega\setminus{\Omega_2}$ and $u_+=\varphi_1$ in $\Omega_1$, where $\varphi_1>0$ denotes an eigenvalue of the Laplace operator in $H^1_0(\Omega_2)$. Choosing $C>0$ big enough, a straightforward argument based on the fact that $0<p<1$ shows that the function $\overline U=Cu_+$ is a super-solution of problem \eq{l31}. Thus, problem \eq{l31} has a nonnegative solution $u$. With the same arguments as in the proof of Lemma \ref{l21} we deduce that $u>0$ in $\Omega$.

\smallskip
{\sc Case 2}: $a<0$.

This case reduces to the previous one. Indeed, by Case 1, let $\underline U$ be the unique solution of the Dirichlet problem
$$\left\{
\begin{array}{lll}
&\di \Delta u-au^p=b(x)f(u)&\qquad\mbox{in $\Omega$}\\
&\di u>0&\qquad\mbox{in $\Omega$}\\
&\di u=\Phi &\qquad\mbox{on $\partial \Omega$.}
\end{array}\right.$$
Then $\underline U$  is a sub-solution of problem \eq{l31}. To construct a super-solution of \eq{l31}, let $G\subset\RR^N$ be an open set with smooth boundary such that $\Omega\subset G$. Let $\varphi_1>0$ be an eigenfunction of the Laplace operator in $H^1_0(G)$. Then $\overline U=M\varphi_1$ is a super-solution of \eq{l31} and $\overline U\geq \underline U$ in $\Omega$, provided that $M>0$ is large enough. The proof of Lemma \ref{l3} is now concluded.
\qed

\medskip
{\bf Proof of Theorem \ref{t1}.}
We first prove the existence of a solution for the nonlinear logistic equation \eq{p} with lower term perturbation. We distinguish two cases, according to the values of the potential function $b(x)$ on $\partial\Omega$. First, if $b>0$ on $\partial\Omega$, then we apply Lemma \ref{l3} for $\Phi\equiv n$. In such a way we obtain an increasing locally bounded sequence of functions that converges to a solution of problem \eq{p}. Next, if $b\geq 0$ on $\partial\Omega$, we apply Lemma \ref{l2} for $\Phi\equiv n$, $h=b+n^{-1}$, $k\equiv a\geq 0$, and $r\equiv 0$.
Now, by Lemma \ref{l2}, we obtain an increasing sequence which is uniformly bounded on every compact subset of $\Omega$. Finally, this sequence converges to a solution of problem \eq{p}. We refer to
C\^{\i}rstea and R\u adulescu \cite[pp.~827-828]{crhouston} for technical details. We also point out that the case $a<0$ can be treated as in the proof of Lemma \ref{l3} by means of a comparison argument.
An alternative argument to establish the existence of a solution of problem \eq{p} if $a<0$ is based on Theorem 1.3 in Dumont, Dupaigne, Goubet and R\u adulescu \cite{ddgr} based on the fact that the mapping $g(x,u):=b(x)f(u)-au$ is a nonnegative smooth function that satisfies the sharpened Keller-Osserman condition.
This concludes the proof of the first part of Theorem \ref{t1}.

\smallskip
Next, we are concerned with the boundary blow-up rate of $u_a$ near $\partial\Omega$. We first observe that relation \eq{63} implies that $h$ is of class $C^2$ in some interval $(0,\delta)$ and $\lim_{t\ri 0^+}h(t)=+\infty$.
We also point out that $h$ is strictly convex near the origin; this follows from
$$\lim_{t\ri 0^+}\frac{h''(t)}{k^2(t)}\, f(h(t)\xi)=\frac{2+\rho\ell_1}{(2+\rho)\xi^{1+\rho}}\qquad\mbox{for all $\xi>0$}.$$
Another direct consequence of this relation is that
$$\lim_{t\ri 0^+}\frac{h(t)}{h''(t)}=\lim_{t\ri 0^+}\frac{h'(t)}{h''(t)}=0\,.$$

Fix  $0<\ep<c$. Our hypotheses imply that there is some $\delta_0>0$ such that $h$ is strictly convex in $(0,\delta_0)$. By continuity, there exists $\delta_1\leq\delta_0$ such that for all $x\in\Omega$ with $d(x)<\delta_1$,
$$\left(c-\ep\right)k^2(d(x))\leq b(x)\leq \left(c+\ep\right)k^2(d(x)).$$
Set
$$\xi^\pm (x):=\left(\frac{2+\rho\ell_1}{(c\mp 2\ep)(2+\rho)} \right)^{1/\rho}.$$
With the same computations as in C\^{\i}rstea and R\u adulescu \cite[pp.~451-452]{crcras} we deduce that
$$\xi^-\leq \liminf_{x\ri x_0}\frac{u_a(x)}{h(d(x))}\leq \limsup_{x\ri x_0}\frac{u_a(x)}{h(d(x))}\leq
\xi^+\,.$$
This implies relation \eq{58}.

At this stage, as soon as we know the blow-up rate of any solution $u_a$ near $\partial\Omega$, it is easy to deduce the uniqueness of the solution. Indeed, let $u$ and $v$ be solutions of problem \eq{p}. Since $\lim_{d(x)\ri 0}u(x)/v(x)=1$, it suffices to apply Lemma \ref{l1} to conclude that $u=v$.
Our proof is now complete.\qed

\medskip
We point out that a stronger existence result than Theorem \ref{t1} holds. More precisely, with the same assumptions as in Theorem \ref{t1}, the nonlinear elliptic problem
$$\left\{\begin{array}{lll}
&\Delta u+au^p+q(x)\,|\nabla u|^\beta=b(x)f(u)&\qquad\mbox{in $\Omega$}\\
&\lim_{x\ri\partial\Omega}u(x)=+\infty\\
&u>0&\qquad\mbox{in $\Omega$}\end{array}\right.$$
has at least one solution, provided that $\beta\in (0,2]$ and $q\in C^{0,\alpha}(\overline\Omega)$ is a non-negative function. The proof combines the arguments from the present paper with those developed by Ghergu and R\u adulescu \cite{grspain}.

\bigskip
{\bf Proof of Theorem \ref{t2}.}
We first observe that for any non-negative function $0\not\equiv \Phi\in
C^{0,\alpha} (\partial \Omega_0)$, the problem
\begin{equation}\label{pp2}\left\{\begin{array}{lll}
&\Delta u+au^p=b(x)f(u)&\qquad\mbox{in $D$}\\
&u=0&\qquad\mbox{on $\partial\Omega$}\\
&u=\Phi&\qquad\mbox{on $\partial \Omega_0$}\\
&u>0&\qquad\mbox{in $D$}\end{array}\right.\end{equation}
has a unique solution. Indeed, let $\overline U$ be the solution of problem \eq{p} if $\Omega$ is replaced by $D$. Then $\overline U$ is a super-solution of problem \eq{pp2} and $\underline U=0$ is a sub-solution, hence \eq{pp2} has at least one solution.
Next, the uniqueness follows by Lemma \ref{l1}.

We now prove that problem \eq{pp} has both a minimal and a maximal solution. Let $u_n$ be the unique solution of problem \eq{pp2} for $\Phi=n$. Thus, by Lemma \ref{l1}, $u_n\leq u_{n+1}$ in $D$. By Theorem \ref{t1}, problem \eq{p} has a solution $u_\infty$ if $\Omega$ is replaced with $D$. Applying again Lemma \ref{l1}, we have $u_n\leq u_\infty$
in $D$. This shows that the sequence $(u_n)$ converges to a solution $\underline u$ of \eq{pp}, which is minimal with respect to other possible solutions.

For all $n\geq 1$ big enough, let $$D_n:=\left\{ x\in D;\ \mbox{dist}\, (x,\partial\Omega_0)>\frac 1n\right\}.$$
Let $w_n$ be the minimal solution of problem \eq{pp} if $D$ is replaced with $D_n$. Thus, by Lemma \ref{l3}, $w_n\geq w_{n+1}$ in $D_n$, which shows that $(w_n)$ converges to to $\overline u$, which is a maximal solution of problem \eq{pp}. A standard regularity arguments that combines Schauder and H\"older estimates ensures that $\underline u$ and $\overline u$ are classical solutions of problem \eq{pp}.

From now on, the proof of Theorem \ref{t2} follows the same lines as in the proof of Theorem \ref{t1}.\qed

\medskip
{\bf Acknowledgments}. The author acknowledges the support by ARRS grants P1-0292-0101 and J1-2057-0101.

\end{document}